\documentclass[12pt]{article}
\usepackage{amsmath,amssymb,amsthm}
\usepackage{mathrsfs}
\usepackage{hyperref}
\usepackage[square, comma, sort&compress, numbers]{natbib}

\allowdisplaybreaks

\makeatletter \@addtoreset{equation}{section}

\makeatletter\@addtoreset{figure}{section}\makeatother

\newtheorem{thm}{Theorem}[section]
\newtheorem{lem}[thm]{Lemma}

\newtheorem{conj}[thm]{Conjecture}
\newtheorem{rem}{Remark}[section]

\begin{document}

\begin{center}
{\large \bf Log-behavior of two sequences related to the elliptic integrals}
\end{center}

\begin{center}
Brian Y. Sun$^{1}$ and James J.Y. Zhao$^{2*}$\\[6pt]

$^{1}$Department of Mathematics and System Science,\\
Xinjiang University, Urumqi 830046, P. R. China\\[6pt]

$^{2}$Center for Applied Mathematics\\
Tianjin University, Tianjin 300072, P. R. China\\[8pt]

Email: $^{1}${\tt brianys1984@126.com},
       $^{2}${\tt jjyzhao@tju.edu.cn}
\\[7pt]

\today
\end{center}

\noindent\textbf{Abstract.}
Two interesting sequences arose in the study of the series expansions of the complete elliptic integrals, which are called the Catalan-Larcombe-French sequence $\{P_n\}_{n\geq 0}$ and the Fennessey-Larcombe-French sequence $\{V_n\}_{n\geq 0}$ respectively. In this paper,
we prove the log-convexity of $\{V_n^2-V_{n-1}V_{n+1}\}_{n\geq 2}$ and $\{n!V_n\}_{n\geq 1}$, the ratio log-concavity of $\{P_n\}_{n\geq 0}$ and the sequence $\{A_n\}_{n\geq 0}$ of Ap\'{e}ry numbers, and the ratio log-convexity of $\{V_n\}_{n\geq 1}$.

\noindent \emph{AMS Classification 2010:} 05A20, 11B37, 11B83


\section{Introduction}

Recently, there is a rising interest in the study of the log-behavior of the
following two sequences defined by
\begin{align}
n^2P_{n}&\!=8(3n^2-3n+1)P_{n-1}-128(n-1)^2P_{n-2}, \label{eq-Pnrec}\\
(n-1)n^2 V_{n}&\!=8(n-1)(3n^2-n-1)V_{n-1}-128(n-2)n^2V_{n-2}, \label{eq-rec}
\end{align}
with the initial values $P_0=V_0=1$ and $P_1=V_1=8$.
The sequences $\{P_n\}_{n\geq 0}$ and $\{V_n\}_{n\geq 0}$ are known as the Catalan-Larcombe-French sequence and the Fennessey-Larcombe-French sequence, respectively. They arise naturally from the series expansions of the complete elliptic integrals, see \cite{Catalan, jv2010, lf2000, lff2002, lff2003}.

Ap\'{e}ry \cite{Apery} introduced the numbers $A_n$, which paly a key role in his proof of the irrationality of $\zeta(3)=\sum_{n=1}^\infty 1/n^3$.
A recurrence relation was also given by Ap\'{e}ry, that is,
\begin{align}\label{eq-An-rec}
n^3 A_n=(34n^3-51n^2+27n-5)A_{n-1}-(n-1)^3 A_{n-2},
\end{align}
with $A_0=1$ and $A_1=5$.

The main objective of this paper is to prove the log-convexity of $\{V_n^2-V_{n-1}V_{n+1}\}_{n\geq 2}$ and $\{n!V_n\}_{n\geq 1}$, the ratio log-concavity of $\{P_n\}_{n\geq 0}$ and $\{A_n\}_{n\geq 0}$, and the ratio log-convexity of $\{V_n\}_{n\geq 1}$. Let us first review some background.

Recall that a real sequence $\{S_n\}_{n \geq 0}$ is said to be log-concave (resp. log-convex) if $S_n^2\geq S_{n-1}S_{n+1}$ (resp. $S_n^2\leq S_{n-1}S_{n+1}$) for all $n\geq 1$, and it is said to be strictly log-concave (resp. strictly log-convex) if the inequality is strict.
Let $\mathcal{L}$ be an operator on $\{S_n\}_{n\geq 0}$ such that
\[ \mathcal{L}(\{S_n\}_{n\geq 0})=\{S_{n-1}S_{n+1}-S_n^2\}_{n\geq 1}.\]
The sequence $\{S_n\}_{n\geq 0}$ is called $k$-log-convex if $\mathcal{L}^i(\{S_n\}_{n\geq 0})$ is log-convex for $0\leq i\leq k-1$, and $\{S_n\}_{n\geq 0}$ is called infinitely log-convex if $\mathcal{L}^k(\{S_n\}_{n\geq 0})$ is log-convex for any $k\geq 1$, see Chen and Xia \cite{ChenXia}.
A real sequence $\{S_n\}_{n\geq 0}$ is called ratio log-concave (resp. ratio log-convex) if the sequence $\{S_{n}/S_{n-1}\}_{n\geq 1}$ is log-concave (resp. log-convex), see Chen, Guo and Wang \cite{CGW2014}.
A real sequence $\{S_n\}_{n\geq 0}$ is called log-balanced if $\{S_n\}_{n\geq 0}$ is log-convex while $\{{S_n}/{n!}\}_{n\geq 0}$ is log-concave, see Do\v{s}li\'{c} \cite{Doslic}.

The log-convexity of $\{P_n\}_{n\geq 0}$, conjectured by Sun \cite{sun2013}, has been proved by Xia and Yao \cite{XY2013} and independently by Zhao \cite{zhao2014}.
The log-concavity of $\{V_n\}_{n\geq 1}$, conjectured by Zhao \cite{zhaocf}, has been confirmed by Yang and Zhao \cite{YangZhao}.
The $2$-log-convexity of $\{P_n\}_{n\geq 0}$ has been shown
by Sun and Wu \cite{SunWu}.
It is natural to consider whether $\{V_n\}_{n\geq 1}$ is $2$-log-concave or not. The first main result of this paper gives the answer.

\begin{thm}\label{th-2lc}
The sequence $\{V_n^2-V_{n-1}V_{n+1}\}_{n\geq 2}$ is strictly log-convex, that is, for $n \geq 3$,
\begin{align}\label{eq-v2lc}
(V_n^2-V_{n-1}V_{n+1})^2 < (V_{n-1}^2-V_{n-2}V_{n})(V_{n+1}^2-V_{n}V_{n+2}).
\end{align}
\end{thm}

Note that Theorem \ref{th-2lc} does not imply the $2$-log-convexity of $\{V_n\}_{n\geq 1}$, since $\{V_n\}_{n\geq 1}$ itself is log-concave.
Chen, Guo and Wang showed that the ratio log-concavity (resp. ratio log-convexity) of a sequence $\{S_n\}_{n\geq N}$ implies the strict log-concavity (resp. strict log-convexity) of the sequence $\{\sqrt[n]{S_n}\}_{n\geq N}$ under an initial condition  \cite[Theorems 3.1 \& 3.6]{CGW2014}.
Although the strictly log-concavity of $\{\sqrt[n]{P_n}\}_{n\geq 1}$ and $\{\sqrt[n]{V_n}\}_{n\geq 1}$ had been proved by Zhao \cite{zhao2016} in a direct way, the ratio log-behaviors of $\{P_n\}_{n\geq 0}$ and $\{V_n\}_{n\geq 1}$ still deserve attention, and are precisely described as follows.

\begin{thm}\label{th-rlcp}
The sequence $\{P_n\}_{n\geq 0}$ is ratio log-concave, that is, for $n\geq 2$,
\begin{align}\label{eq-Prlc}
 (P_n/P_{n-1})^2\geq (P_{n-1}/P_{n-2})(P_{n+1}/P_{n}).
\end{align}
\end{thm}

\begin{thm}\label{th-rlxv}
The sequence $\{V_n\}_{n\geq 1}$ is ratio log-convex, that is, for $n\geq 3$,
\begin{align}\label{eq-Vrlc}
 (V_n/V_{n-1})^2\leq (V_{n-1}/V_{n-2})(V_{n+1}/V_{n}).
\end{align}
\end{thm}

Notice that the strictly log-concavity of $\{\sqrt[n]{P_n}\}_{n\geq 1}$ is a consequence of the criterion \cite[Theorem 3.1]{CGW2014} and Theorem \ref{th-rlcp}, while the strictly log-concavity of $\{\sqrt[n]{V_n}\}_{n\geq 1}$ can not be obtained from the criterion \cite[Theorem 3.1]{CGW2014} and Theorem \ref{th-rlxv}.

Do\v{s}li\'{c} \cite{Doslic} has proved that $\{A_n\}_{n\geq 0}$ is log-convex. The $2$-log-convexity of Ap\'{e}ry numbers has been proved by Chen and Xia \cite{ChenXia}. In this paper, we obtain the ratio log-concavity of the Ap\'{e}ry numbers $A_n$.
\begin{thm}\label{th-rlcA}
The sequence $\{A_n\}_{n\geq 0}$ is ratio log-concave, that is, for $n\geq 2$,
\begin{align}\label{eq-Arlc}
 (A_n/A_{n-1})^2\geq (A_{n-1}/A_{n-2})(A_{n+1}/A_{n}).
\end{align}
\end{thm}
It is easy to check that $\sqrt{A_2}/A_1>\sqrt[3]{A_3}/\sqrt{A_2}$. Thus by Theorem \ref{th-rlcA} and the criterion \cite[Theorem 3.1]{CGW2014}, it follows that
the sequence $\{\sqrt[n]{A_n}\}_{n\geq 1}$ is strictly log-concave, that is, for $n\geq 2$,
\begin{align*}
\left(\sqrt[n]{A_n}\right)^2>\sqrt[n-1]{A_{n-1}}\sqrt[n+1]{A_{n+1}}.
\end{align*}
It should be mentioned that the above inequality was first conjectured by Sun \cite{sun2013}, and then was proved by Luca and St\u{a}nic\u{a} \cite{LucaStanica}.

By further study, we also prove the log-convexity of the sequence $\{n!V_n\}_{n\geq 0}$.

\begin{thm}\label{th-u}
The sequence $\{n!V_n\}_{n\geq0}$ is strictly log-convex, that is, for $n \geq 1$,
\begin{align}\label{eq-nVn}
nV_n^2 < (n+1)V_{n-1}V_{n+1}.
\end{align}
\end{thm}

Since $\{V_n\}_{n\geq 1}$ is log-concave, it follows that the sequence $\{n!V_n\}_{n\geq 1}$ is log-balanced by Theorem \ref{th-u}.
 We notice that Do\v{s}li\'{c}'s criterion of determining log-balancedness \cite[Proposition 3.4]{Doslic} is not available for the sequence $\{n!V_n\}_{n\geq 1}$.
It should be mentioned that Bender and Canfield had given a different criterion \cite[Theorem 1]{BendCanf} for determining log-balancedness of $\{n!S_n\}_{n\geq 1}$, which also does not apply to $\{n!V_n\}_{n\geq 1}$, although they did not name the concept of log-balancedness.

This paper is organized as follows. In Section \ref{s-2}, we prove lower and upper bounds for the ratios $V_n/V_{n-1}$ and $P_n/P_{n-1}$ based on their three-term recurrence relations. The bounds for $A_n/A_{n-1}$, given by Chen and Xia \cite{ChenXia}, are also employed.
These bounds will be used in the proofs of our main results. In Section \ref{s-2lc}, we prove Theorem \ref{th-2lc} by establishing a criterion, which slightly modifies that of Chen and Xia \cite[Theorem 2.1]{ChenXia}.
In Section \ref{S-trr}, we give the proofs of Theorems \ref{th-rlcp}, \ref{th-rlxv} and \ref{th-rlcA} by building two criteria along with the spirit showed in Chen, Guo and Wang \cite[\S 4]{CGW2014}.
In Section \ref{S-u}, we complete the proof of Theorem \ref{th-u}.
We conclude this paper with a few conjectures on log-behaviors related to the Catalan-Larcombe-French sequence and the Fennessey-Larcombe-French sequence.
Since some of the calculations in our proofs are somewhat tedious, we also implement $\mathtt{Maple}$ files to make the checking more convenient.

\section{Bounds for $V_n/V_{n-1}$, $P_n/P_{n-1}$ and $A_n/A_{n-1}$}\label{s-2}

In this section we prove two sets of bounds, one for the ratio $V_n/V_{n-1}$ and the other for the ratio $P_n/P_{n-1}$. A lower bound for $A_n/A_{n-1}$ is also shown. Note that Chen and Xia \cite[\S 4]{ChenXia} have given an upper bound for $A_n/A_{n-1}$, which will be used in our proof. All these bounds are obtained by the heuristic approach shown in \cite[\S 3]{ChenXia} with or without a little polish, and will lead to our main results.
Since three of our main results are related to $V_n$, we first consider the bounds of $V_n/V_{n-1}$. For $n\geq 1$, let
\begin{align}\label{eq-fn}
s(n)=\frac{16(n^5+n^2+3n+12)}{n^5},
\quad {\rm and}\quad
t(n)=\frac{16(n+1)}{n}.
\end{align}

\begin{lem}\label{lem-lub}
Let $s(n)$ and $t(n)$ be given by \eqref{eq-fn}. Then for all integers $n\geq 6$, we have
\begin{align*}
s(n)< \frac{V_n}{V_{n-1}}< t(n).
\end{align*}
\end{lem}

\proof
For notational convenience, let $r(n)=V_n/V_{n-1}$, and
We first prove $r(n)>s(n)$ for $n\geq 6$ by using mathematical induction on $n$. By the recurrence relation \eqref{eq-rec}, we have
\begin{align}\label{eq-gr}
r(n+1)=\frac{8(3n^2+5n+1)}{(n+1)^2}-\frac{128(n-1)}{nr(n)},
\quad n\geq 1,
\end{align}
with the initial value $r(1)=8$.
It is easily checked that $r(6)=20482/1269>1307/81=s(6)$ by \eqref{eq-gr} and \eqref{eq-fn}. Assume $r(n)>s(n)$ holds for $n\geq 6$, and we proceed to show that $r(n+1)>s(n+1)$.
Note that
\begin{align*}
  r(n+1)-s(n+1)
&=\frac{8(3n^2+5n+1)}{(n+1)^2}-\frac{128(n-1)}{nr(n)}\\
&\quad\quad -\frac{16(n^5+5n^4+10n^3+11n^2+10n+17)}{(n+1)^5}\\
&=\frac{8n(n^5+4n^4+5n^3-n^2-12n-33)r(n)-128(n-1)(n+1)^5}{n(n+1)^5r(n)}.
\end{align*}
Clearly, $n^5+4n^4+5n^3-n^2-12n-33=(n^3-1)(n^2+4n+5)-8n-28>0$ for $n\geq 2$  and $r(n)>0$ for $n\geq 1$. By the induction hypothesis, we have $r(n)>s(n)$. Thus for $n\geq 6$, it follows that
\begin{align*}
  r(n+1)-s(n+1)
&>\frac{8n(n^5+4n^4+5n^3-n^2-12n-33)s(n)-128(n-1)(n+1)^5}{n(n+1)^5r(n)}\\
&=\frac{1152(7n^4+5n^3-9n^2-27n-44)}{n^5(n+1)^5r(n)},
\end{align*}
which is clearly positive for $n\geq 6$ since $7n^4+5n^3-9n^2-27n-44=(n+2)(n^2+n+3)(7n-16)+4n^2+11n+52>0$ for $n\geq 6$.
This proves $r(n)>s(n)$ for $n\geq 6$.

For $n\geq 6$, the detailed proof of the inequality $r(n)<t(n)$ are similar to that of $r(n)>s(n)$, and hence is omitted here.
\qed

We now present a lower bound and an upper bound of $P_n/P_{n-1}$. For $n\geq 1$, let
\begin{align}\label{eq-undef}
l(n)=\frac{24 (3 n^2 - 3 n + 1)}{5n^2},\quad {\rm and}\quad
\ell(n)=\frac{16(n^3-n^2-1)}{n^3}.
\end{align}
\begin{lem}\label{lem-rbp}
Let $l(n)$ and $\ell(n)$ be given by \eqref{eq-undef}, then for all integers $n\geq 6$, we have
\begin{align*}
l(n)<\frac{P_n}{P_{n-1}}<\ell(n).
\end{align*}
\end{lem}
\proof
By using mathematical induction on $n$, it is easy to show that $P_n/P_{n-1}>l(n)$ for $n\geq 1$, and $P_n/P_{n-1}<\ell(n)$ for $n\geq 6$. The detailed proof is similar to that of Lemma \ref{lem-lub}, and hence is omitted here.
\qed

In this paper, we adopt the bounds for $A_n/A_{n-1}$ given in \cite{ChenXia}. Let
\begin{align}\label{eq-lbAnn}
p(n)=\frac{34n^3-51n^2+27n-5}{n^3}-\frac{(n-1)^3}{n^3}=\frac{33n^3-48n^2+24n-4}{n^3},
\end{align}
and
\begin{align}\label{eq-ubrAn}
q(n)=17+12\sqrt{2}-\left(\frac{51}{2}+18\sqrt{2}\right)\frac{1}{n}
 +\left(\frac{27}{2}+\frac{609}{64}\sqrt{2}\right)\frac{1}{n^2}
 -\left(\frac{645}{256}+\frac{225}{128}\sqrt{2}\right)\frac{1}{n^3}.
\end{align}

\begin{lem}\label{lem-2.3}
Let $p(n)$ and $q(n)$ be given by \eqref{eq-lbAnn} and \eqref{eq-ubrAn}, respectively. For $n\geq 2$, we have
\begin{align*}
p(n)<\frac{A_n}{A_{n-1}}<q(n).
\end{align*}
\end{lem}

\proof
The inequality $A_n/A_{n-1}<q(n)$ has been proved by Chen and Xia \cite[Lemma 4.1]{ChenXia}.
As noted in \cite[\S 3]{ChenXia}, $p(n)$ is a lower bound for $A_n/A_{n-1}$. It is easy to show that $p(n)<A_n/A_{n-1}$ for $n\geq 2$ by using mathematical induction on $n$. The detailed proof is similar to that of Lemma \ref{lem-lub} and hence is omitted.
\qed

\begin{rem}
Notice that Hou and Zhang \cite{HouZhang} have established an asymptotic method to prove $k$-log-convexity of some sequences except for certain terms at the beginning, and they obtained the bounds by a computer algorithm. With their method, one can obtain the bounds of $S_n/S_{n-1}$ for much more combinatorial sequences $\{S_n\}_{n\geq 0}$.
\end{rem}

\section{Proof of Theorem \ref{th-2lc}}\label{s-2lc}

In this section, we show the proof of Theorem \ref{th-2lc} by presenting a criterion for determining the log-convexity of the sequence $\{S_n^2-S_{n-1}S_{n+1}\}$, where $\{S_n\}_{n\geq 0}$ is a positive sequence that satisfies the recurrence
\begin{equation}\label{eq-trl}
S_n=a(n)S_{n-1}+b(n)S_{n-2}, \quad n\geq 2,
\end{equation}
with real $a(n)$ and $b(n)$.
Our criterion slightly modifies that of Chen and Xia \cite[Theorem 2.1]{ChenXia}.
We notice that in the criterion of Chen and Xia, the sequence $\{S_n\}_{n\geq 0}$ is assumed to be log-convex and an upper bound for $S_n/S_{n-1}$ subject to certain conditions is also needed, while these constrains are not required in ours.

\begin{thm}\label{th-c1}
For a positive sequence $\{S_n\}_{n\geq 0}$ satisfying the relation \eqref{eq-trl}, let
\begin{align*}
c_0(n)=&-b^2(n+1)[a^2(n+2)+b(n+1)-a(n+2)a(n+3)-b(n+3)];\\
c_1(n)=&b(n+1)[2a(n+2)b(n+1)+2a(n+3)a(n+2)a(n+1)\\
&+a(n+3)b(n+2)+2a(n+1)b(n+3)-2a^2(n+2)a(n+1)\\
&-2a(n+2)b(n+2)-3a(n+1)b(n+1)];\\
c_2(n)=&4a(n+1)a(n+2)b(n+1)+2b(n+1)b(n+2)+a^2(n+1)a(n+2)a(n+3)\\
&+a(n+1)a(n+3)b(n+2)+a^2(n+1)b(n+3)-3a^2(n+1)b(n+1)\\
&-a(n+3)a(n+2)b(n+1)-a^2(n+2)a^2(n+1)-b(n+3)b(n+1)\\
&-2a(n+2)a(n+1)b(n+2)-b^2(n+2);\\
c_3(n)=&2a^2(n+1)a(n+2)+2a(n+1)b(n+2)-a(n+1)b(n+3)-a^3(n+1)\\
&-a(n+1)a(n+2)a(n+3)-a(n+3)b(n+2);
\end{align*}
and
\[
\Delta(n)=4c_2^2(n)-12c_1(n)c_3(n).
\]
Suppose that $c_3(n)>0$ and $\Delta(n)>0$ for all $n\geq N$, where
$N$ is a positive integer. If there exists $f(n)$ such that for all
$n\geq N$,
\begin{itemize}
\item[(I)]$\frac{S_n}{S_{n-1}}\geq f(n);$
\item[(II)]$f(n)\geq \frac{-2c_2(n)+\sqrt{\Delta(n)}}{6c_3(n)};$
\item[(III)]$c_3(n)f(n)^3+c_2(n)f(n)^2+c_1(n)f(n)+c_0(n)>0,$
\end{itemize}
then the sequence $\{S_n^2-S_{n-1}S_{n+1}\}_{n\geq N}$ is strictly log-convex, that is, for $n\geq N$,
\begin{equation}\label{eq-lgc2}
(S_{n+1}^2-S_{n}S_{n+2})^2<(S_n^2-S_{n-1}S_{n+1})(S_{n+2}^2-S_{n+1}S_{n+3}).
\end{equation}
\end{thm}

\proof
By the recurrence relation \eqref{eq-trl} and the positivity of the sequence $\{S_n\}_{n\geq 0}$, for $n\geq N$, we have
\begin{align*}
&(S_n^2-S_{n-1}S_{n+1})(S_{n+2}^2-S_{n+1}S_{n+3})-(S_{n+1}^2-S_{n}S_{n+2})^2\\
&\quad=S_{n+1}(2S_nS_{n+1}S_{n+2}+S_{n-1}S_{n+1}S_{n+3}-S_{n+1}^3-S_n^2S_{n+3}-S_{n-1}S_{n+2}^2)\\
&\quad=S_{n+1}(c_3(n)S_n^3+c_2(n)S_n^2S_{n-1}+c_1(n)S_nS_{n-1}^2+c_0(n)S_{n-1}^3)\\
&\quad=S_{n+1}S_{n-1}^3\left[c_3(n)\left(\frac{S_n}{S_{n-1}}\right)^3
+c_2(n)\left(\frac{S_n}{S_{n-1}}\right)^2+c_1(n)\left(\frac{S_n}{S_{n-1}}\right)+c_0(n)\right].
\end{align*}
In order to prove \eqref{eq-lgc2}, it is sufficient to show that for $n\geq N$,
\begin{align}\label{eq-eq1}
c_3(n)\left(\frac{S_n}{S_{n-1}}\right)^3
+c_2(n)\left(\frac{S_n}{S_{n-1}}\right)^2+c_1(n)\left(\frac{S_n}{S_{n-1}}\right)+c_0(n)>0.
\end{align}

Let us consider the polynomial $w(x)=c_3(n)x^3+c_2(n)x^2+c_1(n)x+c_0(n)$. Observe that
$$w'(x)=3c_3(n)x^2+2c_2(n)x+c_1(n).$$
Since $c_3(n)>0$ and $\Delta(n)>0$ for all $n\geq N$, we have the quadratic function $w'(x)\geq 0$
for $x\geq  \frac{-2c_2(n)+\sqrt{\Delta(n)}}{6c_3(n)}$, which
means that $w(x)$ is increasing for $x\in[\frac{-2c_2(n)+\sqrt{\Delta(n)}}{6c_3(n)},+\infty).$
By conditions $(I)$ and $(II)$, we have $\frac{S_n}{S_{n-1}}\geq f(n)\geq\frac{-2c_2(n)+\sqrt{\Delta(n)}}{6c_3(n)}$,
it follows that for $n\geq N$,
\[
w\left(\frac{S_n}{S_{n-1}}\right)\geq w(f(n)).
\]
By condition $(III)$, we have $w(f(n))>0$ for any $n\geq N$. Thus we have $w\left(\frac{S_n}{S_{n-1}}\right)>0$ for $n\geq N$, which leads to \eqref{eq-eq1}.
This completes the proof.
\qed

By the proof of Theorem \ref{th-c1}, one can easily conclude that $\{S_n^2-S_{n-1}S_{n+1}\}_{n\geq N}$ is log-convex if there exists a positive integer $N$ such that for all $n\geq N$, $c_3(n)>0 $, $\Delta(n)<0$, and the conditions $(I)$ and $(III)$ in Theorem \ref{th-c1} holds.

We are now ready to prove Theorem \ref{th-2lc} by using our criterion.

\noindent{\it Proof of Theorem \ref{th-2lc}.}
It is easy to verify that \eqref{eq-v2lc} is true for $n=3,4,5$.
We aim to prove \eqref{eq-v2lc} for $n \geq 6$ by applying Theorem \ref{th-c1}, that is, for $n \geq 6$,
\[
(V_n^2-V_{n-1}V_{n+1})^2 < (V_{n-1}^2-V_{n-2}V_{n})(V_{n+1}^2-V_{n}V_{n+2}).
\]
Compare \eqref{eq-rec} and \eqref{eq-trl}, we have
\[
V_n=a(n)V_{n-1}+b(n)V_{n-2}, \quad n\geq 2,
\]
where
\begin{align*}
a(n)=\frac{8(3n^2-n-1)}{n^2},\quad
b(n)=-\frac{128(n-2)}{n-1}.
\end{align*}

To apply Theorem \ref{th-2lc}, we first verify that $c_3(n)>0$ and $\Delta(n)>0$ for $n\geq 1$.
By computing, it follows that
\begin{equation*}
c_3(n)={\frac {512(n^8+17n^7+131n^6+484n^5+872n^4
+682n^3+51n^2-177n-45)}{(n+1)^6(n+2)^2(n+3)^2}},
\end{equation*}
and
\begin{align*}
\Delta(n)=&\frac{67108864}{(n+1)^8(n+2)^8(n+3)^4 n^2}\left(n^{18}+40n^{17}+752n^{16}+8732n^{15}+69566n^{14}\right.\\
 &+399108n^{13}+1687512n^{12}+5311376n^{11}+12451223n^{10}+21531796n^9\\
&+26834592n^8+23183984n^7+13750782n^6+8285676n^5+10267104n^4\\
&\left.+12477380n^3+9141001n^2+3600576n+596160\right).
\end{align*}
Clearly, both $c_3(n)$ and $\Delta(n)$ are positive for all $n\geq 1$.

Let $N=6$ and $f(n)=s(n)$ for $n\geq N$ where $s(n)$ is defined in \eqref{eq-fn}. We proceed to verify the conditions $(I)$, $(II)$ and $(III)$ in Theorem \ref{th-c1}. It is clear that $\frac{V_n}{V_{n-1}}\geq f(n)$ for $n\geq 6$ by Lemma \ref{lem-lub}, which is the condition $(I)$. We next verify the condition $(II)$. By computing we get
\begin{align*}
&\!\![6c_3(n)f(n)+2c_2(n)]^2-\Delta(n)\\
=&\frac{805306368}{n^{10}(n+3)^4(n+2)^6(n+1)^{12}}
\left(3n^{26}+78n^{25}+952n^{24}+7054n^{23}+37260n^{22}\right.\\
&+172168n^{21}+821087n^{20}+3833124 n^{19}+15316869 n^{18}+49491792 n^{17}\\
&+130518035 n^{16}+295700768 n^{15}+624334735 n^{14}+1306596402 n^{13}\\
&+2645121752 n^{12}+4751027330 n^{11}+6964163254 n^{10}+7754776872 n^9\\
&+5930725839 n^8+2290239180 n^7-689241033 n^6-1426673628 n^5\\
&\left.-697884741 n^4-39615804 n^3+90921852 n^2+32775840 n+3499200\right),
\end{align*}
which is easily checked to be positive for $n\geq 6$.
Note that
\begin{align*}
&6c_3(n)f(n)+2c_2(n)\\
=&\frac{8192}{(n+1)^6(n+2)^4(n+3)^2n^5}\left( n^{15}+22n^{14}+235n^{13}+1362n^{12}+4663n^{11}\right.\\
&+10794n^{10}+23419 n^9+65264 n^8+184207 n^7+395220 n^6+572275 n^5\\
&\left.+497880 n^4+183150 n^3-56592 n^2-67176 n-12960\right ),
\end{align*}
which is clearly positive for $n\geq 6$.
Thus it follows that
\[
6c_3(n)f(n)+2c_2(n)\geq \sqrt{\Delta(n)},
\]
for $n\geq 6$, which is equivalent to the condition $(II)$.

Now it remains to verify the condition $(III)$. To this end, we find that
\begin{align*}
&c_3(n)f(n)^3+c_2(n)f(n)^2+c_1(n)f(n)+c_0(n)\\
=&\frac{3145728}{(n+1)^6(n+2)^4(n+3)^2n^{15}}
\left(66n^{17}+900n^{16}+6674n^{15}+34000n^{14}\right.\\
&+124157n^{13}+336864n^{12}+722550n^{11}+1356276 n^{10}+2548054 n^9\\
&+4990502 n^8+9033247 n^7+13148436 n^6+13877382 n^5+9189072 n^4\\
&\left.+2222712 n^3-1490400 n^2-1178496 n-207360\right)>0
\end{align*}
for all $n\geq 6$, which can be easily checked. This completes the proof.
\qed

\section{Proofs of Theorems \ref{th-rlcp}, \ref{th-rlxv} and \ref{th-rlcA}}\label{S-trr}
In this section we give the detailed proofs of Theorems \ref{th-rlcp}, \ref{th-rlxv} and \ref{th-rlcA}.
Note that Chen, Guo and Wang had established a criterion \cite[Theorem 4.5]{CGW2014} for ratio log-concavity of a sequence subject to the recurrence \eqref{eq-trl}. But their criterion can not be applied to prove our results.
Along with their spirit, we establish two criteria for ratio log-concavity and ratio log-convexity, respectively, of a sequence subject to \eqref{eq-trl}. The first one is as follows.

\begin{thm}\label{thm-crlc}
Let $\{S_n\}_{n\geq 0}$ be a positive sequence satisfying the recurrence relation \eqref{eq-trl},  that is,
\[
S_n=a(n)S_{n-1}+b(n)S_{n-2}, \quad n\geq 2.
\]
Suppose $a(n)>0$ and $b(n)<0$ for $n\geq N$ where $N$ is a nonnegative integer. If there exists two functions $u(n)$ and $v(n)$ such that for all $n\geq N+2$,
\begin{itemize}
\item[$(i)$] $\frac{a(n)}{2}\leq u(n)\leq \frac{S_n}{S_{n-1}}\leq v(n)$;
\item[$(ii)$] $4u^3(n)-3a(n)u^2(n)-a(n+1)b(n)\geq 0$;
\item[$(iii)$] $v^4(n)-a(n)v^3(n)-a(n+1)b(n)v(n)-b(n)b(n+1)\leq0$,
\end{itemize}
then $\{S_n\}_{n\geq N}$ is ratio log-concave, that is, for $n\geq N+2$,
\begin{align}\label{znrlca}
 (S_n/S_{n-1})^2\geq (S_{n-1}/S_{n-2})(S_{n+1}/S_{n}).
\end{align}
\end{thm}

\proof It is clear that \eqref{znrlca} can be rewritten as
\begin{align}\label{znrlc1}
 S_n^3S_{n-2}-S_{n-1}^3S_{n+1}\geq 0.
\end{align}
By the recurrence relation \eqref{eq-trl}, we have
\begin{align*}
  &S_n^3S_{n-2}-S_{n-1}^3S_{n+1}\\
  &\quad =\frac{1}{b(n)}S_n^3(S_n-a(n)S_{n-1})
  -S_{n-1}^3(a(n+1)S_n+b(n+1)S_{n-1})\\
  &\quad= \frac{S_{n-1}^4}{b(n)}
  \left[\left(\frac{S_{n}}{S_{n-1}}\right)^4-
  a(n)\left(\frac{S_{n}}{S_{n-1}}\right)^3
  -a(n+1)b(n)\left(\frac{S_{n}}{S_{n-1}}\right)
  -b(n)b(n+1)\right].
\end{align*}
Note that $b(n)<0$ for $n\geq N+2$. In order to prove \eqref{znrlc1}, it suffices to verify that
\begin{equation}\label{rineq-1}
  \left(\frac{S_{n}}{S_{n-1}}\right)^4-
  a(n)\left(\frac{S_{n}}{S_{n-1}}\right)^3
  -a(n+1)b(n)\left(\frac{S_{n}}{S_{n-1}}\right)
  -b(n)b(n+1)\leq 0,
\end{equation}
for $n\geq N+2$.
Define
  $$h(x)=x^4-a(n)x^3-a(n+1)b(n)x-b(n)b(n+1).$$
Then \eqref{rineq-1} is equivalent to
\[
h\left(\frac{S_n}{S_{n-1}}\right)\leq 0,
\]
for $n\geq N+2$.
Observe that
$$h'(x)=4x^3-3a(n)x^2-a(n+1)b(n),$$
and
$$h''(x)=12x^2-6a(n)x.$$
Since $a(n)>0$, $h''(x)\geq 0$ for $x\geq a(n)/2$, which implies that $h'(x)$ is increasing for $x\geq a(n)/2$.
Note that $u(n)\geq a(n)/2$ by the condition $(i)$. Then we have $h'(x)\geq h'(u(n))$ for $x\geq u(n)$.
By the condition $(ii)$, we have $h'(u(n))\geq 0$. It follows that $h'(x)\geq 0$ for $x\geq u(n)$, and hence $h(x)$ is increasing for $x\geq u(n)$. Then we have $h(S_n/S_{n-1})\leq h(v(n))$ since $u(n)\leq S_n/S_{n-1}\leq v(n)$ by the condition $(i)$.
Now it remains to show that $h(v(n))\leq 0$, which is the condition $(iii)$.
This completes the proof.
\qed

With the help of Theorem \ref{thm-crlc}, we are ready to show the proof of Theorem \ref{th-rlcp}.

\noindent{\it Proof of Theorem \ref{th-rlcp}.}
It is easy to verify that \eqref{eq-Prlc} holds for $2\leq n\leq 5$. We aim to prove \eqref{eq-Prlc} for $n\geq 6$, by applying Theorem \ref{thm-crlc}. Compare \eqref{eq-Pnrec} and \eqref{eq-trl}, we have
\[
P_n=a(n)P_{n-1}+b(n)P_{n-2}
\]
for $n\geq 2$, where
\[
a(n)=\frac{8 (3 n^2 - 3 n + 1)}{n^2},\quad b(n)=-\frac{128 (n - 1)^2}{n^2}.
\]

Set $N=4$. Clearly, $a(n)>0$ and $b(n)<0$ for $n\geq 4$. It suffices to verify the conditions $(i),\ (ii)$ and $(iii)$ in Theorem \ref{thm-crlc}. To this end, let $u(n)=l(n)$ and $v(n)=\ell(n)$ where $l(n)$ and $\ell(n)$ are given by \eqref{eq-undef}. Note that $u(n)=3a(n)/5>a(n)/2$. By Lemma \ref{lem-rbp}, we have $u(n)\leq S_n/S_{n-1}\leq v(n)$ for $n\geq 6$. This verifies the conditions in $(i)$ of Theorem \ref{thm-crlc}. It remains to verify the conditions $(ii)$ and $(iii)$ in Theorem \ref{thm-crlc}.
By computing, we obtain that
\[
4u^3(n)-3a(n)u^2(n)-a(n+1)b(n)=\frac{512\, A(n)}{125 n^6 (n+1)^2},
\]
where
\[
A(n)=21 n^8-21 n^7+229 n^6-1208 n^5+736 n^4+486 n^3-513 n^2+189 n-27,
\]
and
\[
v^4(n)-a(n)v^3(n)-a(n+1)b(n)v(n)-b(n)b(n+1)=\frac{-16384\, B(n)}{n^{12} (n+1)^2},
\]
where
\[
B(n)=4 n^{11}-7 n^{10}-3 n^9-5 n^8+9 n^7+20 n^6+10 n^5-2 n^4-18 n^3-18 n^2-10 n-4.
\]

It is easy to check that $A(n)>0$ and $B(n)>0$ for $n\geq 6$, which confirm the conditions $(ii)$ and $(iii)$ in Theorem \ref{thm-crlc}. This completes the proof.
\qed

As another application of Theorem \ref{thm-crlc}, we now show the proof of Theorem \ref{th-rlcA}.

\noindent{\it Proof of Theorem \ref{th-rlcA}.}
It is easy to verify that \eqref{eq-Arlc} holds for $2\leq n\leq 3$. We aim to prove \eqref{eq-Arlc} for $n\geq 4$, by applying Theorem \ref{thm-crlc}. Compare \eqref{eq-An-rec} and \eqref{eq-trl}, we have
\[
a(n)=\frac{(2n-1)(17n^2-17n+5)}{n^3},\quad b(n)=-\frac{(n-1)^3}{n^3}.
\]

Set $N=2$. Clearly, $a(n)>0$ and $b(n)<0$ for $n\geq 2$. It suffices to verify the conditions $(i),\ (ii)$ and $(iii)$ in Theorem \ref{thm-crlc}. To this end, let $u(n)=p(n)$ and $v(n)=q(n)$ where $p(n)$ and $q(n)$ are given by \eqref{eq-lbAnn} and \eqref{eq-ubrAn}, respectively. Note that $u(n)-a(n)/2=(32n^3-45n^2+21n-3)/(2n^3)$, which is clearly positive for $n\geq 2$. By Lemma \ref{lem-2.3}, we have $u(n)\leq A_n/A_{n-1}\leq v(n)$ for $n\geq 2$. This verifies the conditions in $(i)$ of Theorem \ref{thm-crlc}. It remains to verify the conditions $(ii)$ and $(iii)$ in Theorem \ref{thm-crlc}.
By computing, we obtain that
\[
4u^3(n)-3a(n)u^2(n)-a(n+1)b(n)=\frac{2\, C(n)}{n^9 (n+1)^3},
\]
where
\begin{align*}
C(n)=16352n^{12}-19776n^{11}-29010n^{10}+56240n^9-4659n^8-44808n^7\\
  +31073n^6+1980n^5-11880n^4+6412n^3-1608n^2+192n-8,
\end{align*}
and
\begin{align*}
v^4(n)-a(n)v^3(n)-a(n+1)b(n)v(n)-b(n)b(n+1)=\frac{-3\, D(n)}{4294967296 n^{12} (n+1)^3},
\end{align*}
where
\begin{align*}
 D(n)
\!=&\,(2478196129792+1752346656768\sqrt{2})n^{12}-(6433189920768+4549729320960\sqrt{2})n^{11}\\
& +(4079900164096+2886570344448\sqrt{2})n^{10}+(3923229278208+2773725544448\sqrt{2})n^9\\
& -(7091340886016+5015144103936\sqrt{2})n^8+(3059171226624+2163345012736\sqrt{2})n^7\\
& +(1220059275776+862892127744\sqrt{2})n^6-(1975723880256+1397053488384\sqrt{2})n^5\\
& +(976018184064+690149237616\sqrt{2})n^4-(234159803595+165572939880\sqrt{2})n^3\\
& +(19314604575+13654607400\sqrt{2})n^2+(2591409375+1833597000\sqrt{2})n\\
& -489436875-346275000\sqrt{2}.
\end{align*}

It is easy to check that $C(n)>0$ and $D(n)>0$ for $n\geq 4$, which confirm the conditions $(ii)$ and $(iii)$ in Theorem \ref{thm-crlc}. This completes the proof.
\qed

We now show the criterion for the ratio log-convexity of a sequence  subject to \eqref{eq-trl}.

\begin{thm}\label{thm-crlx}
Let $\{S_n\}_{n\geq 0}$ be a positive sequence satisfying the recurrence relation \eqref{eq-trl}. Suppose $a(n)>0$ and $b(n)<0$ for $n\geq N$ where $N$ is a nonnegative integer. If there exists a function $g(n)$ such that for all $n\geq N+2$,
\begin{itemize}
\item[$(i')$] $\frac{a(n)}{2}\leq g(n)\leq \frac{S_n}{S_{n-1}};$
\item[$(ii')$]  $4g^3(n)-3a(n)g^2(n)-a(n+1)b(n)\geq 0;$
\item[$(iii')$]$g^4(n)-a(n)g^3(n)-a(n+1)b(n)g(n)-b(n)b(n+1)\geq 0,$
\end{itemize}
then $\{S_n\}_{n\geq N}$ is ratio log-convex, that is, for $n\geq N+2$,
\begin{align*}
 (S_n/S_{n-1})^2\leq (S_{n-1}/S_{n-2})(S_{n+1}/S_{n}).
\end{align*}
\end{thm}

\proof
The detailed proof of Theorem \ref{thm-crlx} is similar to that of Theorem \ref{thm-crlc}, and hence is omitted here.
\qed

We are now ready to prove Theorem \ref{th-rlxv}.

\noindent{\it Proof of Theorem \ref{th-rlxv}.\ }
It is easy to check that \eqref{eq-Vrlc} is true for $3\leq n \leq 5$.
We aim to prove \eqref{eq-Vrlc} for $n\geq 6$ by using Theorem \ref{thm-crlx}.
In the proof of Theorem \ref{th-2lc}, we have obtained that
\[
V_n=a(n)V_{n-1}+b(n)V_{n-2},
\]
for $n\geq 2$, where
\[
a(n)=\frac{8(3n^2-n-1)}{n^2},\quad
b(n)=-\frac{128(n-2)}{n-1}.
\]

Let $N=4$. Clearly, $a(n)>0$ and $b(n)<0$ for $n\geq 4$. It suffices to verify the conditions $(i'),\ (ii')$ and $(iii')$ in Theorem \ref{thm-crlx}. For this purpose, let $g(n)=s(n)$ for $n\geq 6$, where $s(n)$ is defined in \eqref{eq-fn}. First by Lemma \ref{lem-lub} we have $g(n)\leq S_n/S_{n-1}$ for $n\geq 6$. Observe that
\[
g(n)-\frac{a(n)}{2}=\frac{4 \left(n^5+n^4+n^3+4 n^2+12 n+48\right)}{n^5}>0,
\]
for $n\geq 1$. This confirms the condition $(i')$ in Theorem \ref{thm-crlx}.

It remains to verify the conditions $(ii')$ and $(iii')$ in Theorem \ref{thm-crlx}. By computation we have that
\[
4g^3(n)-3a(n)g^2(n)-a(n+1)b(n)
 =\frac{1024\, E(n)}{n^{15}(n-1) (n+1)^2},
\]
where
\begin{align*}
E(n)&=n^{18}+3 n^{17}+5 n^{16}+12 n^{15}+48 n^{14}+222 n^{13}+342 n^{12}+300 n^{11}+960 n^{10}\\
 &\quad+2902 n^9+6142 n^8+3956 n^7-448 n^6+9450 n^5+25776 n^4+31536 n^3\\
 &\quad-5184 n^2-48384 n-27648,
\end{align*}
and
\[
g^4(n)-a(n)g^3(n)-a(n+1)b(n)g(n)-b(n)b(n+1)
 =\frac{16384\, F(n)}{ n^{20} (n-1) (n+1)^2},
\]
where
\begin{align*}
F(n)&=24 n^{18}+57 n^{17}+96 n^{16}+234 n^{15}+706 n^{14}+1908 n^{13}+2616 n^{12}+3126 n^{11}\\
 &\quad +8130 n^{10}+18198 n^9+27248 n^8+14970 n^7+5478 n^6+49572 n^5+97308 n^4\\
 &\quad +77760 n^3-58752 n^2-165888 n-82944.
\end{align*}
It is clear that $E(n)> 0$ and $F(n)>0$ for $n\geq 6$. Hence the conditions $(ii')$ and $(iii')$ in Theorem \ref{thm-crlx} are verified
for $n\geq 6$. This completes the proof.
\qed

\section{Proof of Theorem \ref{th-u}}\label{S-u}
In this section, we complete the proof of Theorem \ref{th-u}, the log-convexity of the sequence $\{n!V_n\}_{n\geq 0}$. To make the proof more concise, we need a modified lower bound for the ratio $V_n/V_{n-1}$. For $n\geq 1$, let
\begin{align*}
\tau(n)=\frac{16(n^3+1)}{n^3}.
\end{align*}
Note that $s(n)-\tau(n)=48(n+4)/n^5>0$ for $n\geq 1$ where $s(n)$ is given in \eqref{eq-fn}. Let $r(n)=V_n/V_{n-1}$ and $t(n)$ be defined in \eqref{eq-fn}. Then by Lemma \ref{lem-lub} it is easy to check that
\begin{align}\label{eq-gh}
\tau(n)<r(n)<t(n)
\end{align}
for $n\geq 2$.
With these two bounds, we are now ready to prove Theorem \ref{th-u}.

\noindent{\it Proof of Theorem \ref{th-u}.}
For $n=1$, by the recurrence \eqref{eq-rec}, we have $V_1^2=64<288=2V_0V_1$. We proceed to prove \eqref{eq-nVn} for $n\geq 2$.
Note that \eqref{eq-nVn} can be rewritten as
\[
\frac{r(n)}{r(n+1)}<\frac{n+1}{n}.
\]
Since $r(n)>0$ for $n\geq 1$, by \eqref{eq-gr} we obtain that for $n\geq 1$,
\begin{align*}
\frac{r(n)}{r(n+1)}
 &=\frac{n(n+1)^2r^2(n)}{8n(3n^2+5n+1)r(n)-128(n-1)(n+1)^2}\\
 &=\frac{n(n+1)^2r(n)}{8n(3n^2+5n+1)-128(n-1)(n+1)^2/r(n)},
\end{align*}
with the initial value $r(1)=8$.
Then it suffices to show that
\[
\frac{n(n+1)^2r(n)}{8n(3n^2+5n+1)-128(n-1)(n+1)^2/r(n)}
<\frac{n+1}{n},
\]
for $n\geq 2$.
By \eqref{eq-gh}, we conclude that
\begin{align*}
&\frac{n(n+1)^2r(n)}{8n(3n^2+5n+1)-128(n-1)(n+1)^2/r(n)}
    -\frac{n+1}{n}\\[5pt]
&\;\leq \frac{n(n+1)^2t(n)}{8n(3n^2+5n+1)-128(n-1)(n+1)^2/\tau(n)}
    -\frac{n+1}{n}\\[5pt]
&\;=-\frac{2n^2+n-1}{n(2n^4+2n^3+4n+1)},
\end{align*}
which is clearly negative for $n\geq 2$.
This completes the proof.
\qed

We conclude this paper with a few conjectures related to the Catalan-Larcombe-French sequence and the Fennessey-Larcombe-French sequence.

\begin{conj}\label{conj-1}
The sequence $\{V_n^2-V_{n-1}V_{n+1}\}_{n\geq 2}$ is infinitely log-convex.
\end{conj}
Recently, Wang and Zhu \cite{WangZhu} showed that Stieltjes moment sequences are infinitely log-convex. This provides a possibility for proving Conjecture \ref{conj-1} with some analysis tools.

Let $\mathcal{R}$ be an operator on a sequence $\{S_n\}_{n\geq 0}$ such that
\[ \mathcal{R}(\{S_n\}_{n\geq 0})=\{S_{n+1}/S_{n}\}_{n\geq 0}.\]

\begin{conj}\label{conj-2}
For all integer $k\geq 1$, the sequence $\mathcal{R}^k(\{P_n\}_{n\geq 0})$ except for the first $k$ terms at the beginning is log-concave if $k$ odd, and is log-convex if $k$ even.
\end{conj}

\begin{conj}
For all integer $k\geq 1$, the sequence $\mathcal{R}^k(\{V_n\}_{n\geq 1})$ is log-convex if $k$ odd, and is log-concave if $k$ even.
\end{conj}

\noindent {\bf Acknowledgments.} The authors are  grateful to Arthur L.B. Yang for inspiring discussions and valuable suggestions.
The first author is supported by the young teachers scientific research fund of colleges and universities research program in Xinjiang Uygur Autonomous Region (No. XJEDU2016S032).
This work was supported by the National Science Foundation of China.


{\small
\begin{thebibliography}{10}

\bibitem{Apery} R. Ap\'{e}ry, Irrationalit\'{e} de $\zeta(2)$ et $\zeta(3)$, {\it Ast\'{e}risque} {\bf 61} (1979), 11--13.

\bibitem{BendCanf} E. A. Bender and E. R. Canfield, Log-concavity and related properties of the Cycle Index Polynomials, {\it J. Combin. Theory Ser. A} {\bf 74(1)} (1996), 57--70.

\bibitem{Catalan} E. Catalan, Sur les Nombres de Segner, {\it Rend. Circ. Mat. Palermo} {\bf 1} (1887), 190--201.

\bibitem{CGW2014} W. Y. C. Chen, J. J. F. Guo, and L. X. W. Wang, Infinitely log-monotonic combinatorial sequences, {\it Adv. Appl. Math.} {\bf 52} (2014), 99--120.

\bibitem{ChenXia} W. Y. C. Chen and E. X. W. Xia, The 2-log-convexity of the Ap\'{e}ry Numbers, {\it Proc. Amer. Math. Soc.} {\bf 139} (2011), 391--400.

\bibitem{Doslic} T. Do\v{s}li\'{c}, Log-balanced combinatorial sequences, {\it Int. J. Math. Math. Sci.} {\bf 2005} (2005), 507--522.

\bibitem{HouZhang} Q.-H. Hou and Z.-R. Zhang, Asymptotic $r$-log-convexity and P-recursive sequences, preprint. \href{https://arxiv.org/abs/1609.07840}{arXiv:1609.07840}

\bibitem{jv2010}F. Jarvis and H. A. Verrill, Supercongruences for the Catalan-Larcombe-French numbers, {\it Ramanujan J.} {\bf 22} (2010), 171--186.

\bibitem{lf2000} P. J. Larcombe and D. R. French, On the ``other¡± Catalan numbers: a historical formulation re-examined, {\it Congr. Numer.} {\bf 143} (2000), 33--64.

\bibitem{lff2002} P. J. Larcombe, D. R. French and E. J. Fennessey, The Fennessey-Larcombe-French sequence $\{1, 8, 144, 2432, 40000, ...\}$: formulation and asymptotic form, {\it Congr. Numer.} {\bf 158} (2002), 179-190.

\bibitem{lff2003} P. J. Larcombe, D. R. French and E. J. Fennessey, The Fennessey-Larcombe-French sequence $\{1, 8, 144, 2432, 40000, ...\}$: a recursive formulation and prime factor decomposition, {\it Congr. Numer.} {\bf 160} (2003), 129-137.

\bibitem{LucaStanica} F. Luca and P. St\u{a}nic\u{a}, On some conjectures on the monotonicity of some combinatorial sequences, J. Combin. Number Theory {\bf 4} (2012), 1--10.

\bibitem{SunWu} B. Y. Sun and B. Wu, Two-log-convexity of the Catalan-Larcombe-French sequence, {\it J. Inequal. Appl.} {\bf 2015} (2015), 404.

\bibitem{sun2013} Z. W. Sun, Conjectures involving arithmetical sequences, in Numbers Theory: Arithmetic in Shangri-La, Ser. Number Theory Appl. 8, {\it World Sci., Hackensack, NJ}, 2013, 244--258.

\bibitem{XY2013} E. X. W. Xia and O. X. M. Yao, A criterion for the log-convexity of combinatorial sequences. {\it Electron. J. Combin.} {\bf 20(4)} (2013), \#P3.

\bibitem{WangZhu} Y. Wang and B. X. Zhu, Log-convex and Stieltjes moment sequences, {\it Adv. Appl. Math.} {\bf 81} (2016), 115--127.

\bibitem{YangZhao} A. L. B. Yang and J. J. Y. Zhao, Log-concavity of the Fennessey-Larcombe-French Sequence, {\it Taiwanese J. Math.} {\bf 20(5)} (2016), 993--999.

\bibitem{zhao2014} F.-Z. Zhao, The log-behavior of the Catalan-Larcombe-French sequence, {\it Int. J. Number Theory} {\bf 10} (2014), 177--182.

\bibitem{zhaocf} F.-Z. Zhao, The log-balancedness of combinatorial sequences, The 6th National Conference on Combinatorics and Graph Theory, Guangzhou, November 7--10, 2014, China.

\bibitem{zhao2016} J. J. Y. Zhao, Sun's log-concavity conjecture on the Catalan-Larcombe-French sequence, {\it Acta Math. Sin. (Engl. Ser.)} {\bf 32(5)} (2016), 553--558.

\end{thebibliography}
}
\end{document}